\documentclass{article}
\usepackage{mathrsfs, amsmath, amsthm, amssymb, bm, mathtools, thm-restate}
\usepackage{empheq}
\usepackage{graphicx, xcolor, tikz}
\usetikzlibrary{arrows.meta, intersections, calc, shapes, snakes}
\usepackage{caption, subcaption}
\usepackage{array}
\usepackage{cases}
\makeatletter
\def\th@plain{%
  \upshape 
}
\makeatother



\newtheorem{theorem}{Theorem}[section]
\newtheorem{lemma}[theorem]{Lemma}
\newtheorem{corollary}[theorem]{Corollary}
\newtheorem{conjecture}[theorem]{Conjecture}
\newtheorem*{conjecture*}{Conjecture}
\theoremstyle{definition}

\newtheorem{question}[theorem]{Question}

\usepackage[left=25mm,top=25mm,bottom=25mm,right=25mm]{geometry}
\usepackage[T1]{fontenc}
\usepackage[inline]{enumitem}
\usepackage[square, numbers, sort&compress]{natbib}

\usepackage[pdftex,%
  bookmarks=true,%
  bookmarksnumbered=true, 
  bookmarksopen=true, 
  plainpages=false,%
  pdfpagelabels,%
  colorlinks=true, 
  linkcolor=blue, 
  citecolor=black,%
  anchorcolor=green,
  urlcolor= blue,
  breaklinks=true,
  hyperindex=true]{hyperref}

\newcommand{\etal}{et~al.\ }

\def\mad(#1){\mathrm{mad}(#1)}
\def\int(#1){\mathrm{int}(#1)}
\def\ext(#1){\mathrm{ext}(#1)}
\def\Int(#1){\mathrm{Int}(#1)}
\def\Ext(#1){\mathrm{Ext}(#1)}
\newcommand{\ch}{\mathrm{ch}}

\begin{document}
\title{List star edge coloring of generalized Halin graphs}
\author{Zhengke Miao\textsuperscript{a} \thanks{E-mail address: zkmiao@jsnu.edu.cn}\quad Yimin Song\textsuperscript{b, a} \thanks{E-mail address: ymsong0801@163.com}\quad Tao Wang\textsuperscript{c} \thanks{E-mail address: wangtao@henu.edu.cn}\quad Xiaowei Yu\textsuperscript{a} \thanks{E-mail address: xwyu@jsnu.edu.cn}\\
{\small \textsuperscript{a}Research Institute of Mathematical Science, School of Mathematics and Statistics, }\\
{\small Jiangsu Normal University, Xuzhou, 221116, P. R. China}\\
{\small \textsuperscript{b}School of Mathematical Sciences, }\\
{\small Anhui University, Hefei, 230601, P. R. China}\\
{\small \textsuperscript{c}Institute of Applied Mathematics, }\\
{\small Henan University, Kaifeng, 475004, P. R. China}}
\maketitle
\begin{abstract}
A star $k$-edge coloring is a proper edge coloring such that there are no bichromatic paths or cycles of length four. The smallest integer $k$ such that $G$ admits a star $k$-edge coloring is the star chromatic index of $G$. Deng \etal \cite{MR2933839}, and Bezegov{\'a} \etal \cite{MR3431294} independently proved that the star chromatic index of a tree is at most $\lfloor \frac{3\Delta}{2} \rfloor$, and the bound is sharp. Han \etal \cite{MR3924408} strengthened the result to list version of star chromatic index, and proved that $\lfloor \frac{3\Delta}{2} \rfloor$ is also the sharp upper bound for the list star chromatic index of trees. A generalized Halin graph is a plane graph that consists of a plane embedding of a tree $T$ with $\Delta(T) \geq 3$, and a cycle $C$ connecting all the leaves of the tree such that $C$ is the boundary of the exterior face. In this paper, we prove that if $H \coloneqq T \cup C$ is a generalized Halin graph with $|C| \neq 5$, then its list star chromatic index is at most \[
\max\left\{\left\lfloor\frac{\theta(T) + \Delta(T)}{2}\right\rfloor, 2 \left\lfloor\frac{\Delta(T)}{2}\right\rfloor + 7\right\},
\]
where $\theta(T) = \max_{xy \in E(T)}\{d_{T}(x) + d_{T}(y)\}$. As a consequence, if $H$ is a (generalized) Halin graph with maximum degree $\Delta \geq 13$, then the list star chromatic index is at most $\lfloor \frac{3\Delta}{2} \rfloor$. Moreover, the upper bound for the list star chromatic index is sharp.
\end{abstract}


\section{Introduction}
Unless otherwise stated explicitly, all graphs considered are simple and finite. Given a graph $G$, a \emph{star $k$-edge coloring} of a graph $G$ is a proper edge coloring $\phi: E(G)\rightarrow \{1, 2,\dots, k\}$ such that no path or cycle of length four in $G$ is bichromatic. The \emph{star chromatic index} of $G$, denoted by $\chi'_{\mathrm{star}}(G)$, is the smallest integer $k$ such that $G$ admits a star $k$-edge coloring. Even though star edge coloring attracted much attention, few results are known for general graphs. Let $G$ be a graph with maximum degree $\Delta$. Started from Liu and Deng \cite{MR2416274}, they showed $\chi'_{\mathrm{star}}(G)\le \lceil 16(\Delta-1)^{\frac{3}{2}}\rceil$ when $\Delta\ge 7$ in 2008. In 2013, Dvo\v{r}\'{a}k, Mohar, and \v{S}\'{a}mal \cite{MR3019390} presented a near-linear upper bound $\chi'_{\mathrm{star}}(G)\le  \Delta \cdot 2^{O(1)\sqrt{\log \Delta}}$. For related results, see survey \cite{MR4239912}.

Dvo\v{r}\'{a}k, Mohar, and \v{S}\'{a}mal  pointed out in \cite{MR3019390}: {\em It would be interesting to understand the list version of star edge-coloring}. For a given graph $G$, let $L$ be a list assignment which assigns to each edge $e$ a finite set $L(e)$. If for any list assignment $L$ with $|L(e)| \geq k$ for all $e \in E(G)$, the graph $G$ admits a star edge coloring $\varphi$ such that $\varphi(e)\in L(e)$ for each edge $e \in E(G)$, then $G$ is {\em star $k$-edge choosable}. The \emph{list star chromatic index} of $G$, denoted by $\ch'_{\mathrm{star}}(G)$, is the minimum $k$ such that $G$ is star $k$-edge choosable.

It is difficult to determine the star chromatic index of graphs even for complete graphs and subcubic graphs, just like to determine whether the chromatic index of an arbitrary graph with maximum degree $\Delta$ is $\Delta$ or $\Delta+1$. Lei, Shi, and Song \cite{MR3818598} proved that it is NP-complete to determine whether $\chi'_{\mathrm{star}}(G)\le 3$ for an arbitrary subcubic graph $G$. These motivate scholars to study some special classes of graphs.

In \cite{MR3019390}, the authors studied star edge coloring of subcubic graphs, and showed that if $G$ is a subcubic multigraph, then $\chi'_{\mathrm{star}}(G)\le 7$. Based on this result, they posed the following conjecture and question.

\begin{conjecture} [\cite{MR3019390}]\label{concubic}
If $G$ is a subcubic multigraph, then $\chi'_{\mathrm{star}}(G)\le 6$.
\end{conjecture}

\begin{question}[\cite{MR3019390}]\label{quescubic}
Is it true that $\ch'_{\mathrm{star}}(G)\le 7$ for every subcubic graph $G$? (Perhaps even $\le 6$?)
\end{question}
In 2019, Lu\v{z}ar, Mockov\v{c}iakov\'{a}, and Sot\'{a}k \cite{MR3904837} answered \autoref{quescubic} affirmatively. However, \autoref{concubic} is still open.

\begin{theorem}[\cite{MR3904837}]\label{LISTSUBCUBIC}
If $G$ is a subcubic graph, then $\ch_{\mathrm{star}}'(G) \leq 7$.
\end{theorem}
Lei \etal \cite{MR3818598,MR3764342} studied the subcubic graph with some restrictions on the maximum average degree. Wang \etal \cite{MR3855183} showed that every graph with maximum degree at most $4$ has star chromatic index at most $14$. For results on the list star edge coloring of subcubic graphs can be seen in \cite{MR3826824,MR3802138,MR3763861}.

Two groups of scholars Deng \etal \cite{MR2933839} and Bezegov\'{a} \etal \cite{MR3431294} independently studied the star chromatic index of trees. They proved that for each tree $T$ with maximum degree $\Delta$, we always have $\chi'_{\mathrm{star}}(T)\le \left\lfloor\frac{3\Delta}{2}\right\rfloor$. Moreover, the upper bound is tight. Recently, Han \etal \cite{MR3924408} proved that the upper bound $\left\lfloor\frac{3\Delta}{2}\right\rfloor$ holds for list star chromatic index of trees by ordering edges incident to a common vertex.

\begin{theorem} [\cite{MR3924408}]\label{TREECH}
For each tree $T$ with maximum degree $\Delta$, we always have $\ch'_{\mathrm{star}}(T)\le \left\lfloor\frac{3\Delta}{2}\right\rfloor$. Moreover, this upper bound is tight.
\end{theorem}

The authors in \cite{MR3924408} also presented two upper bounds for list star chromatic index of degenerate graphs by applying some orientation techniques.

The authors in \cite{MR3431294} studied the star edge coloring of outerplanar graphs. Given  an outerplanar graph $G$ with maximum degree $\Delta$, then $\chi'_{\mathrm{star}}(G)\le \left\lfloor\frac{3\Delta}{2}\right\rfloor + 12$ if $\Delta\ge 4$ and $\chi'_{\mathrm{star}}(G)\le 5$ if $\Delta\le 3$. Moreover, the constant 12 in the former case can be decreased to 9 by using more involved analysis. With these two bounds, they proposed the following conjecture.
\begin{conjecture}[\cite{MR3431294}]\label{outerplanar}
Every outerplanar graph $G$ with maximum degree $\Delta\ge 4$ satisfies $\chi'_{\mathrm{star}}(G)\le \left\lfloor\frac{3\Delta}{2}\right\rfloor + 1$.
\end{conjecture}
The best result supporting \autoref{outerplanar} is due to Wang \etal \cite{MR3796353}, they proved the upper bound $\left\lfloor\frac{3\Delta}{2}\right\rfloor+5$ by using an edge-partition technique. In the same paper, the authors studied the star chromatic index of planar graphs and gave the theorem as below.

\begin{theorem}[\cite{MR3796353}]
Let $G$ be a planar graph with maximum degree $\Delta$ and girth $g$. Then
\begin{enumerate}[label = (\alph*)]
  \item $\chi'_{\mathrm{star}}(G)\le \frac{11\Delta}{4} +18$;
  \item $\chi'_{\mathrm{star}}(G)\le \frac{9\Delta}{4} +6$ if $G$ is $K_4$-minor free;
  \item $\chi'_{\mathrm{star}}(G)\le \frac{3\Delta}{2} +18$ if $G$ has no cycle of length four;
  \item $\chi'_{\mathrm{star}}(G)\le \frac{3\Delta}{2} +13$ if $g\ge 5$;
  \item $\chi'_{\mathrm{star}}(G)\le \frac{3\Delta}{2} +3$ if $g\ge 8$.
\end{enumerate}
\end{theorem}

As a special case of planar graphs, Halin graph recently was been widely studied. A \emph{Halin graph} $H \coloneqq T \cup C$ is a plane graph that consists of a plane embedding of a tree $T$ that has no vertices of degree $2$, and a cycle $C$ connecting all the leaves of $T$ such that $C$ is the boundary of the exterior face. The tree $T$ and the cycle $C$ are called the {\em characteristic tree} and the {\em adjoint cycle} of $H$, respectively. If $T$ is relaxed to be a general tree with $\Delta(T) \geq 3$, then $H$ is called a {\em generalized Halin graph}. It is observed that the class of generalized Halin graphs is a superclass of Halin graphs. We call a Halin graph \emph{complete} if all the leaves of $T$ are at the same distance from the root. In the rest of this paper, we always use $H \coloneqq T\cup C$ to denote a generalized Halin graph, in which $T$ and $C$ represent the characteristic tree and the adjoint cycle of $H$, respectively.

As pointed out in \cite{MR3924408}, a strong edge coloring is a special case of star edge coloring. Therefore, the strong chromatic index is an upper bound of the star chromatic index. For the results on strong edge coloring of Halin graphs, we refer the readers to \cite{MR3743945, MR2899886, MR2889505, MR2888094, MR2514906}. Hou \etal \cite{arXiv201014349v1} studied the star chromatic index of complete Halin graphs by showing the following theorem.

\begin{theorem} [\cite{arXiv201014349v1}]\label{halinlarge}
Let $H$ be a complete Halin graph with maximum degree $\Delta\ge 6$. Then
\[
\chi'_{\mathrm{star}}(H)\le \left\lfloor\frac{3\Delta}{2}\right\rfloor+1.
\]
\end{theorem}
Meanwhile, cubic Halin graph was studied by the authors in \cite{arXiv201014349v1} and \cite{MR4242937} independently.

\begin{theorem} [\cite{arXiv201014349v1,MR4242937}]\label{halinsmall}
If $H$ is a cubic Halin graph, then $\chi'_{\mathrm{star}}(H)\le 6$. Furthermore, the bound is tight.
\end{theorem}

Hu \etal \cite{MR3743945} proved that if $H$ is a Halin graph with maximum degree $\Delta \geq 4$, then the strong chromatic index is at most $2\Delta + 1$. In this paper, we prove  the following theorem, which improves \autoref{halinlarge} partially.

\begin{theorem}\label{maintheorem}
If $H$ is a generalized Halin graph with $\Delta \geq 13$, then $\ch_{\mathrm{star}}'(H) \leq \lfloor\frac{3\Delta}{2}\rfloor$. Furthermore, this upper bound is tight.
\end{theorem}
Bezegov\'{a} \etal \cite{MR3431294} pointed out that for every integer $\Delta$, there exists a tree $T$ such that $\chi'_{\mathrm{star}}(T)= \left\lfloor\frac{3\Delta}{2}\right\rfloor$. Let $H \coloneqq T \cup C$ be the generalized Halin graph with the characteristic tree $T$. Under this assumption, we have that $\chi'_{\mathrm{star}}(H)\ge \chi'_{\mathrm{star}}(T)= \left\lfloor\frac{3\Delta}{2}\right\rfloor$. Along with \autoref{maintheorem},  we have that $\ch'_{\mathrm{star}}(H) = \chi'_{\mathrm{star}}(H) = \left\lfloor\frac{3\Delta}{2}\right\rfloor$.

We believe the upper bound in \autoref{maintheorem} does not hold for Halin graphs with small maximum degrees. Although we can not prove the star chromatic index of Halin graph with maximum degree smaller than 13, there is some evidence for this assertion. First, when $\Delta = 3$, according to \autoref{halinsmall}, there exists a cubic Halin graph with star chromatic index 6, which is greater than $\left\lfloor\frac{3 \times 3}{2}\right\rfloor = 4$. Second, when $T$ is a star and $\Delta = 4$, the Halin graph $H$ is the wheel graph $W_{5}$. According to \autoref{wheel}, we have $\chi'_{\mathrm{star}}(H) = 7 > \left\lfloor\frac{3 \times 4}{2}\right\rfloor$.

\begin{theorem} [\cite{Deng2007}]\label{wheel}
Let $W_{n}$ be a wheel graph of order $n \ge 4$. Then
\begin{equation}
\chi'_{\mathrm{star}}(W_{n}) = \left\{\begin{array}{lll}
n+2, &  \mbox{when $n = 5$};\\
n+1, &  \mbox{when $n = 4, 6, 7$};\\
n-1, &  \mbox{when $n \geq 8$}.
\end{array}
\right.\nonumber
\end{equation}
\end{theorem}

To prove \autoref{maintheorem}, we show the following stronger result. The {\em Ore-degree} of an edge $xy$ in a graph $G$ is $\theta(xy) = d(x) + d(y)$. The {\em maximum Ore-degree} of a graph $G$ is defined as $\theta(G) = \max_{xy \in E(G)} \theta(xy) = \max_{xy \in E(G)}\{d(x) + d(y)\}$.
\begin{theorem}\label{MAINRESULT}
If $H \coloneqq T \cup C$ is a generalized Halin graph with $|C| \neq 5$, then \[
\ch_{\mathrm{star}}'(H) \leq \max\left\{\left\lfloor\frac{\theta(T) + \Delta(T)}{2}\right\rfloor, 2 \left\lfloor\frac{\Delta(T)}{2}\right\rfloor + 7\right\}.
\]
\end{theorem}


In order to prove \autoref{MAINRESULT}, we need the concept of distance and the following lemma in the next section. The {\em distance} of two vertices $u$ and $v$ in $G$ is the length of the shortest path connecting $u$ and $v$ in $G$, denoted by $d_{G}(u, v)$.

\begin{lemma}\label{LP}
Let $\Delta$ and $\theta$ be two positive numbers with $\Delta \leq \theta \leq 2\Delta$. Then the optimal value of the following mathematical programming is $\frac{\theta + \Delta}{2}$.
\[
\bgroup
\arraycolsep = 1.4pt
\begin{array}{crccccr}
\mbox{maximize }&\max\{x_{1}&+ &\frac{x_{2}}{2}, &\frac{x_{1}}{2} &+       &x_{2}\}\\
\mbox{subject to }   & x_{1}&+ &x_{2}            &\leq            &\theta  &\\
                     & x_{1}&, &x_{2}            &\leq            &\Delta  &\\
                     & x_{1}&, &x_{2}            &\geq            &0.       &
\end{array}
\egroup
\]
\end{lemma}
In the final section, we slightly improves \autoref{maintheorem} for the complete Halin graphs, and strengthens \autoref{TREECH} for trees in terms of Ore-degree.

\section{Proof of \autoref{MAINRESULT}}\label{Section:2}
Let $H \coloneqq T \cup C$ be a generalized Halin graph, where $T$ is the characteristic tree and $C$ is the adjoint cycle. Let $L$ be an arbitrary list assignment with
\[
|L(e)| \geq k \coloneqq \max\left\{\left\lfloor\frac{\theta(T) + \Delta(T)}{2}\right\rfloor, 2 \left\lfloor\frac{\Delta(T)}{2}\right\rfloor + 7\right\} \mbox{ for every edge $e \in E(H)$}.
\]

Pick a vertex of degree $\Delta$ as the root of $T$, denoted by $r$. If $xy \in E(T)$ and $d_{T}(r, x) = d_{T}(r, y) + 1$, then $y$ is the {\em father} of $x$, and $x$ is a {\em child} of $y$. Note that a vertex may have many children, the root $r$ has no father, and every vertex other than $r$ has exactly one father. We use $\widehat{x}$ to denote the father of $x$. For each vertex $v$, the subgraph induced by $\{uv \mid u \mbox{ is a child of }v\}$ is called a {\em full star} of $v$, denoted by $S_{v}$. Observe that $T$ is a union of some full stars.

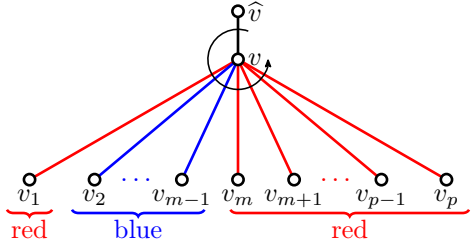
\begin{figure}[h]
\def\s{0.8}
\centering
\begin{tikzpicture}[line width =1pt]
\coordinate (O) at (0, 0);
\coordinate (F) at (90:0.8*\s);
\path [name path = P] (-5*\s, -2*\s)--(5*\s, -2*\s);
\path [name path = OV1] (O)--(210:4*\s);
\path [name path = OVP] (O)--(-30:4*\s);
\path [name path = OV2] (O)--(220:3.2*\s);
\path [name path = OVPM] (O)--(-40:3.2*\s);
\path [name path = OVMM] (O)--(245:2.3*\s);
\path [name path = OVM] (O)--(270:2*\s);
\path [name path = OVMP] (O)--(295:2.3*\s);
\path [name intersections ={of = P and OV1, name = V1, by = {[label = below:$v_{1}$]}}];
\draw[red] (O)--(V1-1);
\path [name intersections ={of = P and OV2, name = V2, by = {[label = below:$v_{2}$]}}];
\draw[blue] (O)--(V2-1);
\path [name intersections ={of = P and OVMM, name = VMM, by = {[label = below:$v_{m-1}$]}}];
\draw[blue] (O)--(VMM-1);
\path [name intersections ={of = P and OVM, name = VM, by = {[label = below:$v_{m}$]}}];
\draw[red] (O)--(VM-1);
\path [name intersections ={of = P and OVMP, name = VMP, by = {[label = below:$v_{m+1}$]}}];
\draw[red] (O)--(VMP-1);
\path [name intersections ={of = P and OVPM, name = VPM, by = {[label = below:$v_{p-1}$]}}];
\draw[red] (O)--(VPM-1);
\path [name intersections ={of = P and OVP, name = VP, by = {[label = below:$v_{p}$]}}];
\draw[red] (O)--(VP-1);
\draw (O)--(F);
\draw[blue, decoration={brace, mirror, raise=.4cm}, decorate] ($(V2-1) - (0.3, 0)$) -- node[below, yshift=-.4cm] {\textcolor{blue}{blue}} ($(VMM-1) + (0.3, 0)$);
\draw[red, decoration={brace, mirror, raise=.4cm}, decorate] ($(VM-1) - (0.1, 0)$) -- node[below, yshift=-.4cm] {\textcolor{red}{red}} ($(VP-1) + (0.3, 0)$);
\draw[red, decoration={brace, mirror, raise=.4cm}, decorate] ($(V1-1) - (0.3, 0)$) -- node[below, yshift=-.4cm] {\textcolor{red}{red}} ($(V1-1) + (0.3, 0)$);
\node[circle, inner sep =1.5, fill=white, draw] () at (V1-1) {};
\node[circle, inner sep =1.5, fill=white, draw] () at (V2-1) {};
\node[circle, inner sep =1.5, fill=white, draw] () at (VMM-1) {};
\node[circle, inner sep =1.5, fill=white, draw] () at (VM-1) {};
\node[circle, inner sep =1.5, fill=white, draw] () at (VMP-1) {};
\node[circle, inner sep =1.5, fill=white, draw] () at (VPM-1) {};
\node[circle, inner sep =1.5, fill=white, draw] () at (VP-1) {};
\node[circle, inner sep =1.5, fill=white, draw] () at (O) {};
\node[circle, inner sep =1.5, fill=white, draw] () at (F) {};
\node[right] at (O) {$v$};
\node[right] at (F) {$\widehat{v}$};
\node () at ($(V2-1)!0.5!(VMM-1)$) {\textcolor{blue}{$\dots$}};
\node () at ($(VMP-1)!0.5!(VPM-1)$) {\textcolor{red}{$\dots$}};
\draw[-{Stealth[length=1.5mm, round]},semithick] (70:0.4) arc (70:360:0.4);
\end{tikzpicture}
\caption{Illustration of $S_{v}$.}
\label{FAMILYSTAR}
\end{figure}

For each vertex $v$ other than $r$, let $v\widehat{v}, vv_{1}, vv_{2}, \dots, vv_{p}$ be all the edges incident with $v$ in the counterclockwise order. Note that when $v$ is a vertex of degree $2$, we have that $p = 1$. We edge-decompose $S_{v}$ into two subgraphs $S_{v, 1}$ and $S_{v, 2}$ such that
\begin{empheq}[left= {E(S_{v, 1}) = \empheqlbrace}]{align*}
&E(S_{v}), &\mbox{when $d(\widehat{v}) + d(v) - 1 - k \leq 0$}, \\
&\{vv_{1}, vv_{m}, \dots, vv_{p}\}, \mbox{where $m = d(\widehat{v}) + d(v) - 1 - k + 2$}, &\mbox{when $d(\widehat{v}) + d(v) - 1 - k \geq 1$},
\end{empheq}
and
\[
E(S_{v, 2}) = E(S_{v}) \setminus E(S_{v, 1}).
\]
We call the edges in $E(S_{v, 1})$ {\em red}, and the edges in $E(S_{v, 2})$ {\em blue}. According to the definition of red edges, it ensures the existence of the red edges $vv_{1}$ and $vv_{p}$ when $E(S_{v, 2}) \neq \emptyset$. An illustration of two types of edges in $S_{v}$, see \autoref{FAMILYSTAR}.

Let $b_{v}$ be the number of blue edges in $S_{v}$. By \autoref{LP}, we have that
\begin{equation}
\textstyle d(\widehat{v}) + (d(v) - 1) - \left\lfloor\frac{d(v)}{2}\right\rfloor \leq k,
\end{equation}
and then
\begin{equation}
\textstyle b_{v} \leq \left\lfloor\frac{d(v)}{2}\right\rfloor.
\end{equation}

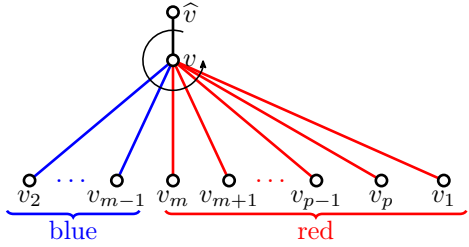
\begin{figure}[h]
\def\s{0.8}
\centering
\begin{tikzpicture}[line width =1pt]
\coordinate (O) at (0, 0);
\coordinate (F) at (90:0.8*\s);
\path [name path = P] (-5*\s, -2*\s)--(5*\s, -2*\s);
\path [name path = OV1] (O)--(4.5*\s, -2*\s);
\path [name path = OVP] (O)--(-30:4*\s);
\path [name path = OV2] (O)--(220:3.2*\s);
\path [name path = OVPM] (O)--(-40:3.2*\s);
\path [name path = OVMM] (O)--(245:2.3*\s);
\path [name path = OVM] (O)--(270:2*\s);
\path [name path = OVMP] (O)--(295:2.3*\s);
\path [name intersections ={of = P and OV1, name = V1, by = {[label = below:$v_{1}$]}}];
\draw[red] (O)--(V1-1);
\path [name intersections ={of = P and OV2, name = V2, by = {[label = below:$v_{2}$]}}];
\draw[blue] (O)--(V2-1);
\path [name intersections ={of = P and OVMM, name = VMM, by = {[label = below:$v_{m-1}$]}}];
\draw[blue] (O)--(VMM-1);
\path [name intersections ={of = P and OVM, name = VM, by = {[label = below:$v_{m}$]}}];
\draw[red] (O)--(VM-1);
\path [name intersections ={of = P and OVMP, name = VMP, by = {[label = below:$v_{m+1}$]}}];
\draw[red] (O)--(VMP-1);
\path [name intersections ={of = P and OVPM, name = VPM, by = {[label = below:$v_{p-1}$]}}];
\draw[red] (O)--(VPM-1);
\path [name intersections ={of = P and OVP, name = VP, by = {[label = below:$v_{p}$]}}];
\draw[red] (O)--(VP-1);
\draw (O)--(F);
\draw[blue, decoration={brace, mirror, raise=.4cm}, decorate] ($(V2-1) - (0.3, 0)$) -- node[below, yshift=-.4cm] {\textcolor{blue}{blue}} ($(VMM-1) + (0.3, 0)$);
\draw[red, decoration={brace, mirror, raise=.4cm}, decorate] ($(VM-1) - (0.1, 0)$) -- node[below, yshift=-.4cm] {\textcolor{red}{red}} ($(V1-1) + (0.3, 0)$);
\node[circle, inner sep =1.5, fill=white, draw] () at (V1-1) {};
\node[circle, inner sep =1.5, fill=white, draw] () at (V2-1) {};
\node[circle, inner sep =1.5, fill=white, draw] () at (VMM-1) {};
\node[circle, inner sep =1.5, fill=white, draw] () at (VM-1) {};
\node[circle, inner sep =1.5, fill=white, draw] () at (VMP-1) {};
\node[circle, inner sep =1.5, fill=white, draw] () at (VPM-1) {};
\node[circle, inner sep =1.5, fill=white, draw] () at (VP-1) {};
\node[circle, inner sep =1.5, fill=white, draw] () at (O) {};
\node[circle, inner sep =1.5, fill=white, draw] () at (F) {};
\node[right] at (O) {$v$};
\node[right] at (F) {$\widehat{v}$};
\node () at ($(V2-1)!0.5!(VMM-1)$) {\textcolor{blue}{$\dots$}};
\node () at ($(VMP-1)!0.5!(VPM-1)$) {\textcolor{red}{$\dots$}};
\draw[-{Stealth[length=1.5mm, round]},semithick] (70:0.4) arc (70:360:0.4);
\end{tikzpicture}
\caption{Illustration of the embedding $\sigma$.}
\label{EMBEDDING}
\end{figure}

We give another embedding $\sigma$ of all the edges incident with $v$, let $v\widehat{v}, vv_{2}, vv_{3}, \dots, vv_{p}, vv_{1}$ be all the edges incident with $v$ in the counterclockwise order, see \autoref{EMBEDDING}. Note that we only change the location of the edge $vv_{1}$ such that all the red edges are consecutive around the vertex $v$. For two distinct edges $e_{1}$ and $e_{2}$ incident with $v$, let
\[
\rho_{v, \sigma}(e_{1}, e_{2}) \coloneqq 1 + |\{e \mid e_{1}, e, e_{2} \mbox{ are in the counterclockwise order around $v$ under the embedding $\sigma$}\}|.
\]

For two distinct edges $e_{1}$ and $e_{2}$ in $E(S_{r})$, let $\sigma$ be the same with the natural embedding of $T$ in $H$, and let $\rho_{v, \sigma}(e_{1}, e_{2})$ be defined as above. By the definition, $\rho_{v, \sigma}(e_{1}, e_{2})$ is different from $\rho_{v, \sigma}(e_{2}, e_{1})$. It is clear that $d(v) = \rho_{v, \sigma}(e_{1}, e_{2}) + \rho_{v, \sigma}(e_{2}, e_{1})$.

Fix an ordering $\preccurlyeq$ of the full stars such that $S_{u} \preccurlyeq S_{w}$ holds only if $d_{T}(r, u) \leq d_{T}(r, w)$. If $S_{u} \preccurlyeq S_{w}$, then we call $S_{u}$ {\em precedes} $S_{w}$. We give an algorithm to color the full stars one by one with the ordering defined by $\preccurlyeq$, and check that no bichromatic $4$-path is produced after each full star is colored.

In the following, we always use $\phi$ to denote the edge coloring in the coloring process, and $\phi(x)$ to denote all the colors on the edges incident with a vertex $x$.

Let $e = v\widehat{v}$ be a pendent edge in $T$. If $e' = u\widehat{u}$ is a pendent edge in another full star $S_{\widehat{u}}$ such that $uv \in E(C)$, then the unordered pair $(e, e')$ is called a {\em bad pair}.
\begin{empheq}[left={\mathcal{B}(e) \coloneqq \empheqlbrace}]{align*}
&\{\phi(e') \mid \mbox{$(e, e')$ is a bad pair}\}, &&\mbox{when $e$ is a pendent edge in $T$};\\
&\emptyset, &&\mbox{otherwise}.
\end{empheq}
From the definition of bad pairs and the structure of generalized Halin graph, we have that $|\mathcal{B}(e)| \leq 2$.

\subsection{List star edge coloring of the characteristic tree}\label{section:2.1}
We first properly $L$-edge color the edges in $S_{r}$. Assume now that $v$ is a vertex other than $r$, and every full star which precedes $S_{v}$ has been colored. We give an algorithm to color the full star $S_{v}$. Let $v\widehat{v}, vv_{1}, vv_{2}, \dots, vv_{p}$ be all the edges incident with $v$ in the counterclockwise order, see \autoref{FAMILYSTAR}. We color the edges in $S_{v}$ in the following order
\[
vv_{1}, vv_{p}, vv_{m}, vv_{m+1}, \dots, vv_{p-1}, vv_{2}, vv_{3}, \dots, vv_{m-1}.
\]
For the red edge $vv_{1}$, we color it with a color from $L(vv_{1}) \setminus (\phi(\widehat{v}) \cup \mathcal{B}(vv_{1}))$. Note that
\[
|L(vv_{1})| - |\phi(\widehat{v})| - |\mathcal{B}(vv_{1})| \geq k - (\Delta + 2) \geq 1,
\]
thus $vv_{1}$ has at least one available color. For the red edge $vv_{p}$, we color it with a color from $L(vv_{p}) \setminus (\phi(\widehat{v}) \cup \mathcal{B}(vv_{p}) \cup \phi(vv_{1}))$. Similarly,
\[
|L(vv_{p})| - |\phi(\widehat{v})| - |\mathcal{B}(vv_{p})| - |\phi(vv_{1})| \geq k - (\Delta + 3) \geq 1,
\]
thus we have at least one available color for $vv_{p}$. According to the definition of red edges, we have that $|\phi(\widehat{v})|$ plus the number of red edges in $S_{v}$ is at most $k$, thus we can properly color all the other red edges in $S_{v}$ such that none of these colors belongs to $\phi(\widehat{v})$. By \autoref{LP}, we have
\begin{equation}
d(v) + \left\lfloor\frac{d(\widehat{v})}{2}\right\rfloor \leq k.
\end{equation}
Thus, for every blue edge $vv_{i}$ in $S_{v}$, we can properly color it with a color in $L(vv_{i})$ but not in
\[
\left\{\phi(e) \mid \mbox{$e$ is incident with $\widehat{v}$ and }\rho_{\widehat{v}, \sigma} (e, v\widehat{v}) \leq \left\lfloor\frac{d(\widehat{v})}{2}\right\rfloor\right\}.
\]

We claim that the two edges in each bad pair have different colors. Suppose to the contrary that a bad pair $(e, e')$ with $\phi(e) = \phi(e')$ is produced after $S_{v}$ is colored. It is observed that one of $e$ and $e'$, say $e = uv$, is contained in $S_{v}$. The edge $e'$ must be in another full star $S_{w}$ which precedes $S_{v}$. By the structure of the generalized Halin graph, $uv$ must be $vv_{1}$ or $vv_{p}$. According to the coloring algorithm, $\phi(uv)$ should be different from $\phi(e')$, a contradiction.

Based on the edge coloring above, the coloring $\phi$ satis?es the following statements.
\begin{enumerate}[label = (\alph*)]
\item For every red edge $uv$ in $S_{v}$, it holds that $\phi(uv) \notin \phi(\widehat{v})$.
\item For every blue edge $uv$ in $S_{v}$, it holds that $\phi(uv) \notin \left\{\phi(e) \mid \mbox{$e$ is incident with $\widehat{v}$ and $\rho_{\widehat{v}, \sigma}(e, v\widehat{v}) \leq \left\lfloor\frac{d(\widehat{v})}{2}\right\rfloor$}\right\}$.
\item The two edges in each bad pair have different colors.
\end{enumerate}

To end this subsection, we show that the resulting edge coloring does not contain bichromatic $4$-paths after the full star $S_{v}$ is colored. Suppose to the contrary that $P$ is a bichromatic $4$-path in $T$. Since we color the full stars one by one according to the fixed ordering defined by $\preccurlyeq$, it is easy to see that $P$ must contain exactly one edge in $S_{v}$. So we may assume that $P = uvxyz$, where $uv \in S_{v}$. It follows that $\phi(uv) = \phi(xy)$ and $\phi(vx) = \phi(yz)$. It is observed that $x$ is the father of $v$. If $uv$ is a red edge, then $\phi(uv) \notin \phi(x)$, this contradicts the fact that $\phi(uv) = \phi(xy)$. Then $uv$ must be a blue edge. (i) Suppose that $y$ is a child of $x$. It follows that $z$ is a child of $y$. Since $\phi(vx) = \phi(yz)$, we have that $yz$ is a blue edge in $S_{y}$ and $\rho_{x, \sigma}(vx, xy) \geq \lfloor\frac{d(x)}{2}\rfloor + 1$. Similarly, $\phi(uv) = \phi(xy)$ implies that $\rho_{x, \sigma}(xy, vx) \geq \lfloor\frac{d(x)}{2}\rfloor + 1$. Hence, $d(x) = \rho_{x, \sigma}(vx, xy) + \rho_{x, \sigma}(xy, vx) \geq 2(\lfloor\frac{d(x)}{2}\rfloor + 1) \geq d(x) + 1$, a contradiction. (ii) Suppose that $y$ is the father of $x$. Since $\phi(vx) = \phi(yz)$, we have that $vx$ is a blue edge in $S_{x}$. Note that there are at most $\lfloor\frac{d(x)}{2}\rfloor$ blue edges in $S_{x}$. It follows that $\rho_{x, \sigma}(xy, vx) \leq b_{x} \leq \lfloor\frac{d(x)}{2}\rfloor$, this implies that $\phi(uv) \neq \phi(xy)$, a contradiction.

Therefore, each full star is colored, no bichromatic $4$-path is produced, and the resulting edge coloring of $T$ is a list star edge coloring.
\subsection{List star edge coloring of the adjoint cycle}
For each edge in $E(C)$, we forbid some colors on specified edges, and guarantee the number of the rest colors is at least three, and then obtain a list star $3$-edge-coloring of $C$.

We define an edge partition $E(C) \coloneqq E_{1} \cup E_{2}$, where
\[
\mbox{$E_{1} \coloneqq \{uv \mid \mbox{$uv \in E(C)$ and $\widehat{u} = \widehat{v}$}\}$, \quad and \quad $E_{2} \coloneqq E(C)\setminus E_{1}$.}
\] Let $uv$ be an arbitrary edge in $E(C)$, define
\begin{empheq}{align*}
\mathscr{A}_{uv}(u) &\coloneqq \mbox{$\{\phi(xy) \mid \mbox{$vuxy$ is a path and $ux \in E(C)$ and $xy \in E(T)$}\}$}, \\
\mathscr{A}_{uv}(v) &\coloneqq \mbox{$\{\phi(xy) \mid \mbox{$uvxy$ is a path and $vx \in E(C)$ and $xy \in E(T)$}\}$}, \\
\mathscr{B}_{u} &\coloneqq \mbox{$\{\phi(e) \mid \mbox{$e$ is incident with $\widehat{u}$ and }\rho_{\widehat{u}, \sigma} (e, u\widehat{u}) \leq \lfloor\frac{d(\widehat{u})}{2}\rfloor\}$}, \\
\mathscr{B}_{v} &\coloneqq \mbox{$\{\phi(e) \mid \mbox{$e$ is incident with $\widehat{v}$ and }\rho_{\widehat{v}, \sigma} (e, v\widehat{v}) \leq \lfloor\frac{d(\widehat{v})}{2}\rfloor\}$}.
\end{empheq}
Since each vertex on $C$ is a pendent vertex of $T$, we have that $|\mathscr{A}_{uv}(u)| = 1$ and $|\mathscr{A}_{uv}(v)| = 1$. It is easy to see that $|\mathscr{B}_{u}| \leq \lfloor\frac{\Delta}{2}\rfloor$, $|\mathscr{B}_{v}| \leq \lfloor\frac{\Delta}{2}\rfloor$. Now we define the available list $L'(uv)$ of $uv$ by
\begin{empheq}[left = {L'(uv) \coloneqq \empheqlbrace}]{align*}
&L(uv) \setminus (\phi(x) \cup \mathscr{A}_{uv}(u) \cup \mathscr{A}_{uv}(v)), &&  \mbox{when $uv \in E_{1}$ and $x = \widehat{u} = \widehat{v}$};\\[0.1cm]
&L(uv) \setminus (\mathscr{A}_{uv}(u) \cup \mathscr{A}_{uv}(v) \cup \mathscr{B}_{u} \cup \mathscr{B}_{v} \cup \phi(u\widehat{u}) \cup \phi(v\widehat{v})), && \mbox{when $uv \in E_{2}$}.
\end{empheq}
It is observed that
\begin{empheq}[left=\empheqlbrace]{align}
    |L'(uv)| &\geq |L(uv)| - \Delta - 2 \geq 3, &\mbox{when $uv \in E_{1}$}; \label{IEQ1}\\
    |L'(uv)| &\geq \textstyle|L(uv)| - 2\left(\left\lfloor\frac{\Delta}{2}\right\rfloor + 2\right) \geq 3, &\mbox{when $uv \in E_{2}$}. \label{IEQ2}
\end{empheq}
Note that $|C| \neq 5$. Wang \etal \cite{MR3855183} proved that $\ch'_{\mathrm{star}}(C) = 3$ if $C$ is a cycle of length $|C| \neq 5$. Therefore, the adjoint cycle $C$ has a star $L'$-edge coloring.

\subsection{The resulting coloring is a list star edge coloring of $H$}
Suppose to the contrary that $uvxyz$ be a bichromatic $4$-path or $4$-cycle ($u = z$ in this case) in $H$. That is, $\phi(uv) = \phi(xy)$ and $\phi(vx) = \phi(yz)$.

By the above discussion, we may assume that $\{uv, vx, xy, yz\} \cap E(T) \neq \emptyset$ and $\{uv, vx, xy, yz\} \cap E(C) \neq \emptyset$.

(i) {\bf Assume that $\bm{x \notin V(C)}$}. Then $vx$ and $xy$ are edges in $T$. By symmetry, we may assume that $uv \in E(C)$. It is observed that $vx$ is a pendent edge in $T$ and $x = \widehat{v}$. If $uv$ is an edge in $E_{1}$, then $\phi(uv) \notin \phi(x)$, $\phi(uv) \neq \phi(xy)$, a contradiction. So we may assume that $uv$ is an edge in $E_{2}$. Since $\phi(xy) = \phi(uv) \notin \mathscr{B}_{v}$, we have that
\begin{equation}\label{xy-xv}
\rho_{x, \sigma}(xy, xv) \geq \left\lfloor\frac{d(x)}{2}\right\rfloor + 1,
\end{equation}
and hence
\begin{equation}\label{xv-xy}
\rho_{x, \sigma}(xv, xy) \leq d(x) - \left\lfloor\frac{d(x)}{2}\right\rfloor - 1 \leq \left\lfloor\frac{d(x)}{2}\right\rfloor.
\end{equation}
Since $\phi(vx) = \phi(yz)$, we have that $yz$ cannot be in $E(C)$. This implies that $yz \in E(T)$.

Assume that $y$ is a child of $x$. Since $\phi(vx) = \phi(yz)$, we have that $yz$ is a blue edge in $S_{y}$ and $\rho_{x, \sigma}(xv, xy) \geq \lfloor\frac{d(x)}{2}\rfloor + 1$, this contradicts the inequality \eqref{xv-xy}. So we may assume that $y$ is the father of $x$. According to the coloring algorithm for $T$, $\phi(vx) = \phi(yz)$ implies that $vx$ is a blue edge in $S_{x}$. Note that there are at most $\lfloor\frac{d(x)}{2}\rfloor$ blue edges in $S_{x}$, so we have $\rho_{x, \sigma}(xy, xv) \leq b_{x} \leq \lfloor\frac{d(x)}{2}\rfloor$, this contradicts the inequality \eqref{xy-xv}.

(ii) {\bf Assume that $\bm{x \in V(C)}$}. Suppose that $vx$ and $xy$ are two edges in $E(C)$. Then one of $uv$ and $yz$, say $uv$, is in $E(T)$. It is observed that $\phi(xy) = \phi(uv) \in \mathscr{A}_{xy}(x)$, a contradiction. Since $x$ is a vertex of degree three in $H$, we may assume that one of $vx$ and $xy$ is in $E(C)$, and the other one is in $E(T)$. Without loss of generality, let $vx \in E(C)$ and $xy \in E(T)$. If $uv \in E(C)$, then $\phi(uv) = \phi(xy) \in \mathscr{A}_{uv}(v)$, a contradiction. It follows that $uv \in E(T)$. Since $\phi(uv) = \phi(xy)$, we have that $u \neq y$. Then $uv$ and $xy$ are in a bad pair, and they should have different colors, this contradicts that $\phi(uv) = \phi(xy)$.

Therefore, $\phi$ is a star $L$-edge coloring of $H$.

\section{Final discussions}
The class of complete Halin graphs is a subclass of generalized Halin graphs. Let $H \coloneqq T \cup C$ be a complete Halin graph. Following the discussions in \autoref{Section:2}, each component of $H[E_{2}]$ must be an isolated edge. Let $uv$ be an arbitrary edge in $E_{2}$. It is observed that either $\mathscr{A}_{uv}(u) \subseteq \mathscr{B}_{u}$ or $\mathscr{A}_{uv}(v) \subseteq \mathscr{B}_{v}$, so the inequality \eqref{IEQ2} can be replaced with $|L'(uv)| \geq |L(uv)| - 2\lfloor\frac{\Delta}{2}\rfloor - 3 \geq 3$, that requires $|L(uv)| \geq 2\lfloor\frac{\Delta}{2}\rfloor + 6$. Then we have the following result for complete Halin graphs.

\begin{theorem}
If $H \coloneqq T \cup C$ is a complete Halin graph with $|C| \neq 5$, then \[
\ch_{\mathrm{star}}'(H) \leq \max\left\{\left\lfloor\frac{\theta(T) + \Delta(T)}{2}\right\rfloor, 2 \left\lfloor\frac{\Delta(T)}{2}\right\rfloor + 6\right\}.
\]
\end{theorem}
\begin{corollary}
If $H$ is a complete Halin graph with $\Delta \geq 11$, then $\ch_{\mathrm{star}}'(H) \leq \left\lfloor\frac{3\Delta}{2}\right\rfloor$.
\end{corollary}

Faudree \etal \cite{MR1412876} proved that the strong chromatic index of a tree $T$ is $\theta(T) - 1$. When $T$ is a path, a special tree, it is easy to check that the list star chromatic index does not exceed $\left\lfloor\frac{\theta(T) + \Delta(T)}{2}\right\rfloor$. When $T$ is a tree with $\Delta(T) \geq 3$, following the arguments in \autoref{section:2.1} by changing only a few words (the bad pairs are not required), we can prove the next result in terms of maximum Ore-degree.

\begin{theorem}\label{ORE}
If $T$ is a tree, then $\ch'_{\mathrm{star}}(T) \leq \left\lfloor\frac{\theta(T) + \Delta(T)}{2}\right\rfloor$.
\end{theorem}

When there are two adjacent vertices of degree $\Delta(T)$, the maximum Ore-degree is $2\Delta(T)$, thus $\ch'_{\mathrm{star}}(T) \leq \left\lfloor\frac{2\Delta(T) + \Delta(T)}{2}\right\rfloor = \left\lfloor\frac{3\Delta(T)}{2}\right\rfloor$, which is consistent with the bound in \autoref{TREECH}. Thus, \autoref{ORE} strengthens \autoref{TREECH}. Observe that the upper bound in \autoref{ORE} is tight.

\vskip 3mm \vspace{0.3cm} \noindent{\bf Acknowledgments.} This work was supported by NSFC grants: 11971205, 12031018, 11901252.


\end{document}